\documentclass{article}
\usepackage{amsmath, amsthm, amsfonts, amssymb}
\usepackage{hyperref}
\newtheorem{theorem}{Theorem}[section]

\newtheorem{proposition}[theorem]{Proposition}

\theoremstyle{definition}
\newtheorem{definition}[theorem]{Definition}
\newtheorem{example}[theorem]{Example}

\theoremstyle{remark}

\numberwithin{equation}{section}

\begin{document}
\title{\bf Maximal Ideals in a Bicomplex Algebra and Bicomplex Gelfand-Mazur Theorem}
\date{\textbf{ Kulbir Singh and Romesh Kumar }}
\vspace{0in}
\maketitle

\textbf{Abstract:} In this paper we study the maximal ideals in a commutative ring of bicomplex numbers and then we describe the maximal ideals in a bicomplex algebra. We found that the kernel of a nonzero multiplicative $\mathbb {BC}$-linear functional in a commutative bicomplex Banach algebra need not be a maximal ideal. Finally, we introduce the notion of bicomplex division algebra and generalize the Gelfand-Mazur theorem for the bicomplex division Banach algebra.\\\\
$\textbf{Keywords:}$ \;Bicomplex modules, hyperbolic-valued norm, maximal ideals, bicomplex algebra.\\\\ 

\begin{section} {Inroduction and Preliminaries}
Bicomplex numbers are being studied for a long time. In the last several years, the theory of bicomplex numbers has enjoyed a renewed interest and bicomplex functional analysis have been studied in \cite{YY}, \cite{HH}, \cite{G_2G_2}, \cite{Hahn} and references therein. Due to the presence of zero divisors in the ring of bicomplex numbers, the theory of bicomplex analysis differs in many aspects from analysis of one complex variable. For instance, the spectrum of a bounded operator on a Banach $\mathbb {BC}$-module is unbounded (see, \cite [Corollary 3.11] {HH}).\\\\
In this section we summarize some basic properties of bicomplex numbers. The ring of bicomplex numbers is defined as follows:
\begin{align}\mathbb{B}\mathbb{C}&=\left\{\;Z=z+jw\;|\;z,\;w\in \mathbb{C}(i)\;\right\},
\end{align} where  $i$ and $j$ are two imaginary units such that
$ ij=ji, \;i^2=j^2=-1$ and $\mathbb{C}(i)$ is the set of complex numbers with the imaginary unit $i$.   
 In (1.1), if $z=x$ is real and $w=iy$ is a purely imaginary number with $ij=k$, then we obtain the ring of hyperbolic numbers 
 $$\mathbb{D}=\left\{x+ky:k^2=1 \;\text{and}\;x,\;y\in \mathbb{R}\; \text{with}\;k \notin \mathbb R\right\}.$$
 Since bicomplex numbers are defined as the pair of two complex numbers connected through another imaginary unit, there are several natural notions of conjugations. Let $Z=z+jw\in \mathbb{B}\mathbb{C}$. Then the following three conjugates can be defined in $\mathbb{B}\mathbb{C}$:\\
(i)\;$\overline {Z}=\overline{z}+j\overline{w}$\;,\;\; (ii)\;$Z^{\dagger}=z-jw$\;,\;\;(iii)\;$Z^{*}= \overline{z}-j\overline{w},$\;\;
where $\overline{z}$, $\overline{w}$ denote the usual complex conjugates to $z, w$ in $\mathbb C(i).$  
With each kind of the three conjugations, three possible moduli of a bicomplex number are defined as follows:\\ 
(i) $|Z|^2_j=Z\;.\;\overline {Z},$
(ii) $|Z|^2_i=Z\;.\;Z^{\dagger},$  
(iii) $|Z|^2_k=Z\;.\;Z^{*}.$\\\\ 	  The system of bicomplex numbers $\mathbb {BC}$ is considerably simplified by the introduction of two numbers $e_1=\frac{1+ij}{2} \;\;\text{ and its $\dagger$-conjugate}\;\; e_2= e_1^\dagger =\frac{1-ij}{2}\;.$
The numbers $e_1$ and $e_2$ are zero divisors in $\mathbb {BC}$. We denote the set of all zero divisors in $\mathbb{B}\mathbb{C}$ by $\mathcal {NC}$ read as null cone and is defined as $$\mathcal {NC} =\left\{Z=z+jw \in \mathbb {BC}\;|\; Z\neq 0,\; z^2+w^2=0 \right\}.$$ In fact, $e_1$ and $e_2$ are hyperbolic numbers and satisfy the following properties:
$$e_1^2=e_1,\; e_2^2=e_2,\; e_1^*=e_1,\;e_2^*=e_2,\;e_1+e_2=1,\;e_1.\;e_2=0.$$ The numbers $e_1$ and $e_2$ are called the idempotent basis of $\mathbb {BC}$. Thus, every bicomplex number $Z=z+jw$ can be written in a unique way as:  
\begin{align} \label{iddec} Z&=e_1z_1+e_2z_2 \;,
\end{align}  where $z_1=z-iw$ and $z_2=z+iw$ are elements of  $\mathbb{C}(i)$. Formula $(\ref{iddec})$ is called the idempotent representation of a bicomplex number $Z$. 
Further, the two sets $e_1\mathbb {BC}$ and $e_2\mathbb {BC}$ are (principal) ideals in the ring $\mathbb {BC}$ such that $e_1\mathbb {BC} \; \cap \;e_2\mathbb {BC}= \left\{0\right\}$. Thus, we can write $\mathbb {BC}= e_1\mathbb {BC}+e_2\mathbb {BC}$ and is called the idempotent decomposition of $\mathbb {BC}$. The hyperbolic-valued or $\mathbb D$-valued norm $|Z|_k$ of a bicomplex number $Z=e_1z_1+e_2z_2$ is defined as $|Z|_k=e_1|z_1|+e_2|z_2|.$ \\	 
A $\mathbb {BC}$-module $\mathcal A$ can be written as 
$\mathcal A=e_1\mathcal A_1+e_2\mathcal A_2$ and is called the idempotent decomposition of $\mathcal A$, where $\mathcal A_1=e_1\mathcal A$ and $\mathcal A_2=e_2\mathcal A$ are $\mathbb C(i)$-vector spaces.  
 Assume that $\mathcal A_1$ and $\mathcal A_2$ are normed spaces with respective norms $\|.\|_1$ and $\|.\|_2$. For any $x \in \mathcal A$, set
\begin{align} \|x\|_{\mathbb D}=\|e_1x_1+e_2x_2\|_{\mathbb D}=e_1\|x\|_1+e_2\|x_2\|_2 \;.\end{align} Then $\|.\|_{\mathbb D}$ defines a hyperbolic-valued or $\mathbb D$-valued norm on $\mathbb{B}\mathbb{C}$-module $\mathcal A$.\\ For the following Definitions, one can refer to \cite{KS2} and \cite{GKZ}.
\begin{definition} A bicomplex algebra $\mathcal A$ is a module over $\mathbb {BC}$  with a multiplication defined on it which satisfy the following conditions:\\ $(i)$ $x(y+z)=xy+xz$\; 
$(ii)$ $(x+y)z=xy+yz$\;
$(iii)$ $x(yz)=(xy)z$\\ 
$(iv)$ $\lambda(xy)=(\lambda x)y=x(\lambda y)$ \;for all $x,\;y,\;z \in \mathcal A$ and all scalars  $\lambda\in \mathbb {BC}$.
\end{definition}

\begin{definition} A bicomplex algebra $\mathcal A$ with a $\mathbb D$-valued norm $\|.\|_\mathbb D$ relative to which $\mathcal A$ is a Banach module such that for every  $x,y$ in $\mathcal A$,   
$$ \|xy\|_\mathbb D \leq' \|x\|_\mathbb D\|y\|_\mathbb D$$ 
is called a $\mathbb D$-normed bicomplex Banach algebra.
\end{definition}
\begin{definition} Let $\mathcal A$ be a bicomplex algebra. Then a bicomplex homomorphism $f:\mathcal A\rightarrow \mathbb {BC}$ is a $\mathbb {BC}$-linear functional such that $$f(xy)=f(x)f(y), \;\;\forall\; x,y \in \mathcal A.$$
\end{definition} 
For more details on bicomplex analysis and various properties of bicomplex Banach algebras, one can refer to \cite{YY}, \cite{G_1G_1}, \cite{G_2G_2}, \cite{KS2}, \cite{Hahn},  \cite{GKZ}, \cite{KK}, \cite{RR}, \cite{XX}, and many other references therein.
\end{section}

\section  {Maximal Ideals in a Bicomplex Algebra}
 
In this section, we define the bicomplex ring homomorphisms in the ring $\mathbb {BC}$ and then we find the maximal ideals in $\mathbb {BC}$ as well as in a bicomplex algebra. It is important to know the behaviour of the bicomplex homomorphisms when they are evaluated on invertible elements, as the behaviour of a bicomplex homorphism is not same as that of a complex homomorphism. For example, if $\phi$ is a complex homomorphism on a complex algebra $\mathcal A$ with identity $u$, then $\phi(u)=1$ and $\phi(x)\neq 0$ for every invertible $x\in \mathcal A$. Now for bicomplex homomorphisms, M. E. Luna-Elizarraras et al. \cite{GKZ} proved the following result  :
\begin{theorem} If $\phi$ is a bicomplex homomorphism on a bicomplex algebra $\mathcal A$ with identity $u$ and $\phi$ is such that at least one of its values is an invertible bicomplex number, then $\phi(u)=1$ and for every invertible $x\in \mathcal A$ its image $\phi(x)$ is invertible with $\phi(x^{-1})=(\phi(x))^{-1}.$ 
\end{theorem}
Now here we first define bicomplex ring homomorphisms and then we will show $e_1\mathbb {BC}$ and $e_2\mathbb {BC}$ are two maximal ideals in $\mathbb {BC}$. Consider the ring $\mathbb {BC}$. Let $f:\mathbb {BC} \rightarrow \mathbb {BC}$ be a map such that $f(x+y)=f(x) +f(y)\;\;\text{and}\;\;f(xy)=f(x)f(y), \;\;\forall\; x,y \in \mathcal A,$ then $f$ is called a bicomplex ring homomorphism.\\\\
It is well known that $1,0, e_1$ and $e_2$ are only four idempotent elements in the ring $\mathbb {BC}$. Further, if $f:\mathbb {BC}\rightarrow \mathbb {BC}$ is a bicomplex ring homomorphism then an idempotent element will be mapped to an idempotent element. Thus, the possible bicomplex ring homomorphisms in the ring $\mathbb {BC}$ are defined as follows:\\\\
 (i) $f(Z)= Z$, $\forall \;Z\in \mathbb {BC}$, i.e., the identity homomorphism. In this case $Ker f=\left\{0\right\}$, the zero ideal in $\mathbb {BC}$. Moreover, if $f(Z)=\overline{Z}$ or $Z^{\dagger}$ or $Z^{*}$, then $f$ is also a ring homomorphism with $Ker f=\left\{0\right\}$.\\\\
(ii) $f(Z)=0$, $\forall \;Z\in \mathbb {BC}$, i.e., the zero homomorphism. Here $Ker f= \mathbb {BC}$ which is an improper ideal in $\mathbb {BC}$.\\\\
(iii) Define $f(Z)=e_1Z$ ( or $e_2Z$), $\forall \;Z\in \mathbb {BC}$. Then $f$ is a bicomplex ring homomorphism with $Ker f= e_2\mathbb {BC}$ ( or $e_1\mathbb {BC}$), which is a proper ideal in $\mathbb {BC}$.\\\\
Thus, $\mathbb {BC}$ and $\left\{0\right\}$ are two trivial ideals and $e_1\mathbb {BC}$ and $e_2\mathbb {BC}$ are two proper ideals in the ring $\mathbb {BC}.$ Since  $\left\{0\right\}$ is contained in both $e_1\mathbb {BC}$ and $e_2\mathbb {BC}$, the two proper ideals $e_1\mathbb {BC}$ and $e_2\mathbb {BC}$ are maximal ideals can be proved either by showing $\mathbb {BC}/e_l\mathbb {BC}$, $(l=1,2)$, is a field or by using the fundamental theorem of ring homomorphisms.
\begin{proposition} In the ring $\mathbb {BC}$ of bicomplex numbers, the ideals $e_1\mathbb {BC}$ and $e_2\mathbb {BC}$ are maximal.
\end{proposition}
\begin{proof}
Let $I_1=e_1\mathbb {BC}$ and $I_2=e_2\mathbb {BC}$ be two proper ideals in the ring $\mathbb {BC}$. Consider the quotient set 
 \begin{align*}  \mathbb {BC}/I_1&=\left\{ I_1+Z\;|\;Z\in \mathbb {BC}\right\}=\left\{ I_1+e_2z_2\;|\;Z=e_1z_1+e_2z_2\in \mathbb {BC}\right\}.
 \end{align*}
 For any $Z=e_1z_1+e_2z_2, W=e_1w_1+e_2w_2\in \mathbb {BC}$, we define addition and multiplication on $\mathbb {BC}/I_1$ as follows:
 $$(I_1+Z) + (I_1+W)= I_1+(Z+W)= I_1+e_2(z_2+w_2)\;\;\text{and}$$
 $$(I_1+Z)(I_1+W)= I_1+(ZW)= I_1+e_2(z_2w_2).$$
Then it is easy to show that $\mathbb {BC}/I_1$ is a commutative ring with multiplicative identity $I_1+e_2$. Now for any $Z=e_1z_1+e_2z_2\in \mathbb {BC}$ with $z_2 \neq 0$, $I_1+e_2z_2$ is a nonzero element in $\mathbb {BC}/I_1$. Further, $z_2\in \mathbb C(i)$, there exists $0\neq w\in \mathbb C(i)$ such that $z_2w=wz_2=1$. This implies that $I_1+e_2w$ is a nonzero element in $\mathbb {BC}/I_1$ with $$(I_1+e_2w) (I_1+e_2z_2)=I_1+e_2=(I_1+e_2z_2)(I_1+e_2w).$$ Thus, every nonzero element in $\mathbb {BC}/I_1$ has a multiplicative inverse. Hence $\mathbb {BC}/I_1$ is a field which implies that $I_1$ is a maximal ideal in the ring $\mathbb {BC}$. In the same manner, we can show that $I_2$ is a maximal ideal in the ring $\mathbb {BC}$.\\\\  
\textbf{Alternative Proof} : 
 Define a map $f:\mathbb {BC}\rightarrow \mathbb C(i)$ as $$f(Z)=f(z+jw)=z+iw,\;\; \forall \;Z\in \mathbb {BC}.$$ Clearly, $f$ is an onto ring homomorphism.\\ Then, by fundamental theorem of ring homomorphisms, 
\begin{align} \label{ftrh} \mathbb {BC}/Ker f \cong \mathbb C(i).\end{align}
 Now for any $Z\in \mathbb {BC}$,\begin{align*} Z=z+jw\in Ker f&\Leftrightarrow f(Z)=f(z+jw)=0\\
  &\Leftrightarrow z+iw=0\\
 & \Leftrightarrow Z=z+jiz \Leftrightarrow Z\in I_1
\end{align*}
implies that $Ker f= I_1$ and thus, by using (\ref{ftrh}),  $\mathbb {BC}/I_1$ is a field. Hence $I_1$ is a maximal ideal in the ring $\mathbb {BC}$.\\ In the same way, on defining an onto ring homomorphism $g:\mathbb {BC}\rightarrow \mathbb C(i)$ by $$g(Z)=g(z+jw)=z-iw,$$ we find that $Ker\;g=I_2$ is a maximal ideal in the ring $\mathbb {BC}$.
\end{proof}
\begin{theorem}\label{MI} Let $\mathcal A$  be a  bicomplex algebra and $\mathcal A=e_1\mathcal A_1+e_2\mathcal A_2$ be its idempotent decomposition. Then every maximal ideal in $\mathcal A$ is either of the form \begin{align} \label{mid1} e_1\mathcal A_1+e_2 J_2,\;\text{where}\; J_2 \;\text{is a maximal ideal in}\; \mathcal A_2\ \end{align} $~~~~~~~~~~~~~~~~~~~~~~~~~~~~~~~~~~~~~~~~~~~or$ \begin{align} \label{mid2} e_1J_1+e_2\mathcal A_2, \;\text{where}\; J_1 \;\text{is a maximal ideal in}\; \mathcal A_1.\end{align} 
\end{theorem}
\begin{proof} Suppose $J=e_1\mathcal A_1+e_2 J_2$ is not a maximal ideal in $\mathcal A$, where $J_2$ is maximal in $\mathcal A_2$. Then $J$ will contained in some proper ideal $J'$ of $\mathcal A$. Clearly, $J'$ will be of the form $e_1\mathcal A_1+e_2 J'_2$, where $J'_2$ is an ideal in $\mathcal A_2$ and $J'_2$ contains $J_2$ which is a contradiction to the maximality of $J_2$. Hence ideals of the form (\ref{mid1}) are maximal in $\mathcal A$. Similarly one can show that ideals defined in (\ref{mid2}) are maximal in $\mathcal A$. Further, let us suppose that $J$ be any other maximal ideal in $\mathcal A$. Then by the idempotent decomposition of $J$, we can write $J=e_1J_1+e_2J_2$. If $J_1= \left\{0\right\}$ or $J_2= \left\{0\right\}$, then in this case $J$ cannot be maximal as it is properly contained either in $A_1$ or in $A_2$. Thus, we assume that both $J_1$ and $J_2$ are nonzero ideals. Then, again in this case, it is clear that $J$ belong to the category (\ref{mid1}) or to (\ref{mid2}). Hence every maximal ideal in a bicomplex algebra $\mathcal A$ is of the form (\ref{mid1}) or of (\ref{mid2}).            
\end{proof}
 In a commutative complex Banach algebra $\mathcal A$ with identity, kernel of every nonzero complex homomorphism is a maximal ideal in $\mathcal A$. Further, every maximal ideal in $\mathcal A$ is a kernel of some nonzero complex homomorphism. Thus, there exists a one-to-one correspondence between nonzero complex homomorphisms and maximal ideals in a complex Banach algebra.\\\\  
\textbf{Question}: Does the above correspondence between nonzero $\mathbb {BC}$-homomorphisms and maximal ideals exist for a Banach $\mathbb {BC}$-algebra?\\
Here we provide an example to show that this type of correspondence may not exists in a Banach $\mathbb {BC}$-algebra.\\
\begin{example} Consider a commutative Banach $\mathbb {BC}$-algebra $\mathbb {BC}$ with identity $e_1+e_2=1$. We now define a map $f:\mathbb {BC} \rightarrow \mathbb {BC}$ as $f(Z)=Z$. Then $f$ is a $\mathbb {BC}$-homomorphism with $ker f=\left\{0\right\}$, i.e., the zero ideal. But $\left\{0\right\}$ is not a maximal ideal in $\mathbb {BC}$ as it is properly contained in $e_1\mathbb {BC}$ and $e_2\mathbb {BC}$.
 Thus, the correspondence between nonzero $\mathbb {BC}$-homomorphisms and maximal ideals in a Banach $\mathbb {BC}$-algebra is not one-to-one.\end{example}  

\section {Bicomplex Division Algebra}
Here in this section, we introduce the notion of bicomplex division algebra and extend the famous Gelfand-Mazur Theorem for bicomplex division algebras. 
  \begin{definition} An algebra $\mathcal A$ over $\mathbb {BC}$ is said to be a bicomplex division algebra if every nonzero element which does not belong to the null cone has a multiplicative inverse in $\mathcal A$. That is, every nonzero element in $\mathcal A -\mathcal{NC}$ is invertible.  
  \end{definition}
 \begin{example} The ring $\mathbb {BC}$ forms a bicomplex division algebra.
 \end{example}
 In \cite{GKZ}, it is proved that if $\mathcal A$ is a bicomplex Banach algebra, then the sets $e_1 \mathcal A$ and $e_2 \mathcal A$ are complex Banach algebras. Further, it is shown that if $z$ is an invertible element in $\mathcal A$ with identity $u$, then $e_1z$ is invertible in $e_1 \mathcal A$ with identity $e_1u$ and $e_2z$ is invertible in $e_2 \mathcal A$ with identity $e_2u$. Thus, we have the following results: 
 
\begin{proposition} \label{bcdvnal7390}Let $\mathcal A$ be a bicomplex division algebra with identity $u$. Then the sets $e_1 \mathcal A$ and $e_2 \mathcal A$ are complex division algebras.
\end{proposition}
\begin{proof} Let $z$ be a nonzero element in  $e_1 \mathcal A$. Choose some nonzero element $w \in e_2 \mathcal A$ so that $Z=e_1z+e_2w \in \mathcal A$. Then $Z$ being a nonzero element in $\mathcal A - \mathcal{NC}$ is invertible. Thus, there exists a nonzero element $Z'=e_1z'+e_2w' \in \mathcal A-\mathcal{NC}$ such that $ZZ'=u=Z'Z$. Clearly, $z' \in e_1 \mathcal A$ with $zz'=e_1u=z'z$. Thus, every nonzero element in $e_1 \mathcal A$ is invertible. Hence $e_1 \mathcal A$ is a complex division algebra. In the same way, one can show that$e_2 \mathcal A$ is a complex division algebra.
\end{proof} 
\begin{proposition} Let $\mathcal A$ be a bicomplex algebra such that the sets $e_1 \mathcal A$ and $e_2 \mathcal A$ are complex division algebras. Then $\mathcal A$ is a bicomplex division algebra.
\end{proposition}
\begin{proof} Let $Z$ be a nonzero element in $\mathcal A-\mathcal{NC}$. Then we can write $Z=e_1z_1+e_2z_2$, where $z_1=e_1Z\in e_1 \mathcal A$ and $z_2=e_2Z \in e_2 \mathcal A$ with $z_1, z_2 \neq 0$. Since $ e_1 \mathcal A$ is a complex division algebra, there exists a nonzero element $z_1' \in  e_1 \mathcal A$ such that $z_1z_1'=e_1u=z_1'z_1$. Similarly, there exists  some $0 \neq z_2' \in  e_2 \mathcal A$ such that $z_2z_2'=e_2u=z_2'z_2$. Clearly, $Z'=e_1z_1'+e_2z_2'$ is a nonzero element in $\mathcal A-\mathcal{NC}$ with $ZZ'=u=Z'Z$. Thus, every nonzero element in $\mathcal A - \mathcal{NC}$ is invertible implies that $\mathcal A$ is a bicomplex division algebra.    
\end{proof}

\begin{definition} Let $\mathcal A$ be a bicomplex Banach algebra with identity $u$. Then for any $x\in \mathcal A$, the spectrum of $x$ is defined as $$\sigma_{(x)}=\left\{ \lambda \in \mathbb {BC} \;:\; x-\lambda u \text{ \;is\; not\; invertible\; in\; $\mathcal A$}\right\}.$$ 
\end{definition}
The spectrum of a bounded $\mathbb {BC}$-linear operator on a Banach $\mathbb {BC}$-module is studied by F. Colombo, I. Sabadini and D. C. Struppa in \cite{HH} and they proved the following result: 
\begin{theorem} \label{splop}  Let $\mathcal B(V)$ be a Banach $\mathbb {BC}$-module of all bounded operators acting on $V$. Then for any $T\in \mathcal B(V)$, the spectrum $\sigma (T)$ is unbounded.
\end{theorem}
 We consider the Theorem \ref{splop} in a more general setting and obtained the following result:
 \begin{theorem} Let $\mathcal A$ be a bicomplex division algebra with identity $u$. Then for any $x \in \mathcal A$, the spectrum $\sigma_{(x)}$ of $x$ is unbounded.
\end{theorem}
\begin{proof}  By Proposition \ref{bcdvnal7390}, for any $x\in \mathcal A$, we can write $x=e_1x_1+e_2x_2$, where $x_l$ belongs to the complex division algebra $e_l \mathcal A$, $l=1,2$. Then
\begin{align*} \sigma_{(x)}&=\left\{ \lambda \in \mathbb {BC} \;:\; x-\lambda u \text{ is\; not\; invertible\; in}\; \mathcal A\right\}\\
&= \left\{ \lambda \in \mathbb {BC} \;:\; x-\lambda u\; \in \;\;\mathcal {NC} \cup \left\{0 \right\}\; \right\}\\
&=\left\{ e_1\lambda_1+e_2 \lambda_2 \in \mathbb {BC} : e_1(x_1-\lambda_1 u_1)+e_2(x_2-\lambda_2 u_2) \in \mathcal {NC} \cup \left\{0 \right\} \right\}
\end{align*}
which implies that atleast one between $ x_1-\lambda_1 u_1$  and  $x_2-\lambda_2 u_2$ is a zero vector. In the first case, if $ x_1-\lambda_1 u_1=0$, then $x_1=\lambda_1 u_1$, while $\lambda_2 \in \mathbb C(i)$ as $x_2-\lambda_2 u_2$ can be any vector in $\mathbb C(i)$. Thus, in this case, we obtain $$\sigma_{(x)}= e_1 \lambda_1 + e_2 \mathbb C(i)=e_1 \sigma_{(x_1)} + e_2 \mathbb C(i).$$ Similarly, if $ x_2-\lambda_2 u_2=0$, then $x_1-\lambda_1 u_1$ can be any vector in $\mathbb C(i)$ implies that $x_2=\lambda_2 u_2$ while $\lambda_1 \in \mathbb C(i)$. Thus, in this case, we obtain $$\sigma_{(x)}= e_1 \mathbb C(i)+ e_2 \lambda_2= e_1 \mathbb C(i)+ e_2 \sigma_{(x_2)}.$$ 
Thus, in both cases, the spectrum of $x$ is given by $$\sigma_{(x)}= (e_1 \sigma_{(x_1)} + e_2 \mathbb C(i)) \cup (e_1 \mathbb C(i)+ e_2 \sigma_{(x_2)}).$$  This shows that $\sigma_{(x)}$ is unbounded. 
\end{proof}For a complex Banach division algebra $\mathcal A$, the spectrum of every element $x$ in $\mathcal A$ contains only one element. Thus, by defining a map $x\rightarrow \sigma_{(x)}$, we see that $\mathcal A$ is isometrically isomorphic to $\mathbb C$ which is the famous Gelfand-Mazur Theorem.
 We now generalize the Gelfand-Mazur Theorem for bicomplex division algebras.

\begin{theorem} Let $\mathcal A$ be a $\mathbb D$-normed bicomplex Banach division algebra. Then $\mathcal A$ is $\mathbb D$-isometrically isomorphic to $\mathbb {BC}$.
\end{theorem}
\begin{proof} Let $\mathcal A$ be a $\mathbb D$-normed bicomplex Banach division algebra. Then $e_1\mathcal A$ and $e_2\mathcal A$ are complex division Banach algebras with norms $\|.\|_1$ and $\|.\|_2$ respectively induced by the $\mathbb D$-valued norm on $\mathcal A$. Further, for any $x=e_1x_1+e_2x_2\in \mathcal A$, the spectrum $\sigma_{(x)}$ of $x$ in $\mathcal A$ is given by $$\sigma_{(x)}= (e_1 \sigma_{(x_1)}+ e_2 \mathbb C(i)) \cup (e_1 \mathbb C(i)+e_2 \sigma_{(x_2)}),$$ where $\sigma_{(x_l)}$ denote the spectrum of $x_l$ in $e_l \mathcal A$, $l=1,2$. Also, for each $x_l \in e_l\mathcal A$, $\sigma_{(x_l)}$ consists of exactly one point say $\lambda_l$, $l=1,2.$ We now define a map $f:\mathcal A\rightarrow \mathbb {BC}$ as $$f(x)= e_1 \sigma_{(x_1)}+ e_2 \sigma_{(x_2)}, \;\forall \;x \in \mathcal A.$$ That is, $f(e_1x_1+e_2x_2)= e_1 \lambda_1+e_2 \lambda_2$. Then it is easy to show that $f$ is a $\mathbb {BC}$-isomorphism of $\mathcal A$ onto $\mathbb {BC}$. Further, for each $x\in \mathcal A$, \begin{align*}  
|f(x)|_k&=|\lambda|_k= \|\lambda u\|_\mathbb D=\|e_1 \lambda_1 u_1+e_2 \lambda_2 u_2\|_\mathbb D \\
&=e_1\|\lambda_1u_1\|_1+e_2 \|\lambda_2 u_2\|_2\\
&=e_1\|x_1\|_1+e_2\|x_2\|_2=\|x\|_\mathbb D \end{align*}
Hence $f$ is a $\mathbb D$-isometry. 
\end{proof}

\noindent Kulbir Singh, \textit{Department of Mathematics, Govt. Degree College Kathua, Jammu,   J\&K - 184101, India.}\\
E-mail :\textit{ singhkulbir03@yahoo.com}\\

\noindent Romesh Kumar, \textit{Department of Mathematics, University of Jammu, Jammu, J\&K - 180006, India.}\\
E-mail :\textit{ romesh\_jammu@yahoo.com}\\

 \end{document}